\documentclass[11pt]{amsart}
\usepackage{amsmath}
\usepackage{amssymb}
\usepackage{mathtools}
\usepackage{tikz}

\usepackage{amsmath}
\usepackage{amssymb}
\usepackage{mathtools}
\usepackage{tikz}

\newcommand{\mb}{\mathbb}
\newcommand{\be}{\begin{equation}}
\newcommand{\ee}{ \end{equation}}

\newcommand{\dd}{\,\mathrm{d}}

\newcommand{\vn}[1]{\|#1\|}
\newcommand{\vm}[1]{\left|#1\right|}

\newcommand{\Cross}{$\mathbin{\tikz [x=1.4ex,y=1.4ex,line width=.2ex, red] \draw (0,0) -- (1,1) (0,1) -- (1,0);}$}%

\newcommand{\pt}{\partial}
\newcommand{\xx}{\frac{n-1}{2}}

\newcommand{\ls}{\lesssim}
\newcommand{\norm}[1]{\left\lVert#1\right\rVert}
\newcommand{\defeq}{\vcentcolon=}

\newcommand{\lmd}{\lambda}

\newtheorem{theorem}{Theorem}[section]

\theoremstyle{definition}

\theoremstyle{remark}
\newtheorem{remark}[theorem]{Remark}

\numberwithin{equation}{section}

\begin{document}

\title[]{A deterministic counterexample for high dimensional $ L^2 L^{\infty} $ Strichartz estimates for the Wave equation}

\author{Cristian Gavrus}

\curraddr{}
\email{gcd2006@gmail.com}
\thanks{}


\date{}

\begin{abstract}
In this note we discuss the question of 
homogeneous $ L^2 L^{\infty} $ Strichartz estimates for the Wave equation in dimensions $ n \geq 4 $ raised by Fang and Wang and recently shown to fail by 
Guo, Li, Nakanishi and Yan using probability theory. 
We record a deterministic example for disproving this estimate. 
\end{abstract}

\maketitle

\section{Introduction}

The question of determining the validity of the following statement 
was raised in \cite[Remark 1]{fang2006some} and was recently answered negatively in \cite{guo2018boundary}
\be \label{W} \tag{W}
 \vn{e^{it \vm{D}} u}_{L^2_t L^{\infty}_x} \lesssim  \vn{u}_{\dot{H}^{\frac{n-1}{2}}} \quad \text{holds for all } u \in \mathcal{S}(\mathbb{R}^n) \text{ when }  n \geq 4. 
\ee 
The function $ v = e^{it \vm{D}} u $ solves the wave equation $ (- \pt_t^2 + \Delta ) v= 0 $. Here  $ \mathcal{S}(\mathbb{R}^n) $ denotes the space of Schwartz functions, while $ \dot{H}^s $ denote homogeneous Sobolev spaces. 

In fact \cite{guo2018boundary} disproved estimates $ e^{it \vm{D}^a} : \dot{H}^{\frac{n-a}{2}} \to L^2_t L^{\infty}_x $ 
for general operators given by $ a \in (0,2] $, using stable Levy processes. The case $ a=2 $ corresponds to the Schr\"{o}dinger equation.

In this note we record a proof of the failure of \eqref{W} using basic Fourier analysis. This argument is specific to the wave equation and does not apply to the Schr\"{o}dinger case.

\

\begin{theorem} \label{ThmW}
The statement \eqref{W} is not true. 
\end{theorem}

\

Estimates of $ L^q_t L^r_x $ norms of solutions to dispersive equations with initial conditions in $ \dot{H}^s $ spaces are known as Strichartz estimates, see \cite{strichartz1977restrictions}, \cite{yajima1987existence}, \cite{kapitanski1989some}, \cite{ginibre1992smoothing}, \cite{lindblad1995existence}, \cite{ginibre1995generalized},  \cite{sogge1995lectures}, \cite{tao2006nonlinear}.
The exponents $ (q,r) = (2,\bar{r}) $ when $ r $ is finite and minimal  
\footnote{i.e. $ \bar{r} = \frac{2(n-1)}{n-3} $ for the Wave equation ($ n \geq 4 $) and  $ \bar{r} = \frac{2n}{n-2} $ for the Schr\"{o}dinger eq. ($ n \geq 3 $). } (blue point in Figure 1 Left)
are called the endpoints, and have been famously resolved in \cite{keel1998endpoint}. We refer to \cite[Section 2.3, see "Knapp example"]{tao2006nonlinear} for a discussion of the minimality of $ \bar{r} $. 

The case $ r= \infty $ is usually omitted when stating Strichartz estimates in higher dimensions ($ n \geq 4 $ for the Wave equation and $ n \geq 3 $ for the Schr\"{o}dinger equation), see for example \cite{tao2006nonlinear}, \cite{klainermanlecture}, \cite{sogge1995lectures}, because the strongest and most useful cases occur on the blue line in Fig. 1, i.e. the range between  $ L_t^2 L_x^{\bar{r}} $ and $ L_t^{\infty} L_x^2 $. Other cases are typically obtained by Sobolev embeddings. The $ L_t^2 L_x^{\infty} $ exponents differ from the endpoints, unless $ n=3 $ for the Wave equation or $ n=2 $ for the Schr\"{o}dinger equation, Figure 1 Right. In those low-dimension cases, counterexamples are given in \cite{klainerman1993space}, \cite{montgomery1998time}. See Remark \ref{Wn3} below. The present case is the red point in Fig. 1 Left. 

The case $ r= \infty $ is generally delicate because $ \dot{W}^{\frac{n}{r},r} \nsubseteq L^{\infty} $, i.e. Sobolev embeddings fail into $ L^{\infty} $. Thus \eqref{W}  does  not follow from the Keel-Tao endpoint $ L_t^2 L_x^{\bar{r}} $ in  \cite{keel1998endpoint}, and instead one can only conclude from  \cite{keel1998endpoint}  estimates that are true 
if one replaces
 $ L^{\infty}_x $ with $ BMO $ or Besov spaces (or $ L_t^2 L_x^{\infty} $ estimates restricted to frequency localized functions) - see the ending of Remark \ref{Wn3}. 

\begin{figure} \label{fig1} 


\begin{tikzpicture}[scale=1.8] 

\draw [step=0.5,thin,gray!40] (0,0) grid (2.4,2.4);

\draw [->] (0,0) -- (2.5,0) node [below] {$  \frac{1}{r}$};
\draw [->] (0,0) -- (0,2.5) node [left] {$ \frac{1}{q} $};

\draw [thin,black] (0,1.5) -- (0.8,1.5) -- (1.5 ,0) -- (0,0) -- cycle;

\draw [thick,blue] (0.8,1.5) -- (1.5 ,0);

\fill [blue] (0.8,1.5) circle (1.5 pt);
\fill [black] (1.5 ,0) circle (1.5 pt);
\fill [red] (0,1.5) circle (1.5 pt);

\draw [black] (1.5,1.5)node [above] {$ L^2 L^{\bar{r}} $ enpoint };
\draw [black] (1.5 ,0) node [below] {$ \frac{1}{2} $};
\draw [black] (1.8 ,0) node [above] {$ L^{\infty} L^{2} $};
\draw [black] (0,1.5) node [left] {$ \frac{1}{2} $};
\draw [black] (0.3,1.6)  node [above] {$ L^{2} L^{\infty}  $};
\draw [black] (0,0.7) node [left] {$ L^{q} L^{\infty}  $};
\node   {\Cross};
\end{tikzpicture}
\qquad \qquad
\begin{tikzpicture}[scale=1.8] 
\draw [step=0.5,thin,gray!40] (0,0) grid (2.4,2.4);

\draw [->] (0,0) -- (2.5,0) node [below] {$  \frac{1}{r}$};
\draw [->] (0,0) -- (0,2.5) node [left] {$ \frac{1}{q} $};

\draw [thin,black] (0,1.5) -- (1.5 ,0) -- (0,0) -- cycle;
\draw [thick,blue] (0,1.5) -- (1.5 ,0);

\fill [black] (1.5 ,0) circle (1.5 pt);
\fill [blue](0,1.5) circle( 1.5 pt);
\fill [red] (0,1.5) circle ( 0.7 pt);

\draw [black] (1.5 ,0) node [below] {$ \frac{1}{2} $};
\draw [black] (0,1.5) node [left] {$ \frac{1}{2} $};
\draw [black] (0.3,1.6)  node [above] {$ L^{2} L^{\infty}  $};
\draw [black] (1.8 ,0) node [above] {$ L^{\infty} L^{2} $};
\node   {\Cross};

\end{tikzpicture}
\caption{Left: high dimensions.   \qquad   Right: $ n=3 $ for W.E.} \label{fig}
\end{figure}
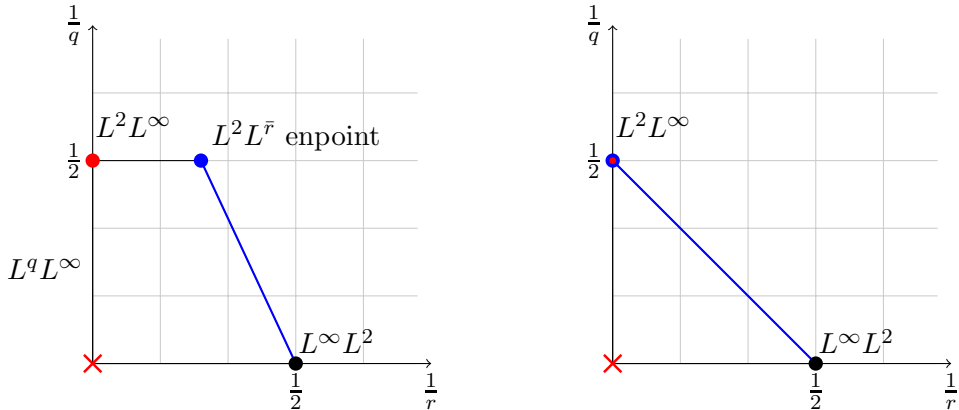

Estimates $ e^{i t \vm{D}} : \dot{H}^{\frac{n}{2}- \frac{1}{q}} \to L_t^{q} L_x^{\infty} $ for $ q \in (2,\infty) $ 
are proved in \cite{klainerman1998algebraic}, \cite{fang2006some}. The case $ q=\infty $ trivially fails at $ t=0$ due to $ \dot{H}^{\frac{n}{2}}  \nsubseteq L^{\infty}_x $. The case that remained was $ q=2 $, which was recently treated in \cite{guo2018boundary}.  

\begin{remark} \label{Wn3}
When $ n=3 $ the estimate in \eqref{W} is disproved in \cite[Remark 1]{klainerman1993space} (also presented in \cite{klainermanlecture}), exhibiting  concentration along null rays. Using the $ 3d $ fundamental solution (i.e. Kirchhoff's formula) they show there exists a sequence of solutions to $ (- \pt_t^2 + \Delta ) \phi_n = 0 $
with normalized initial velocities $ \vn{ \pt_t \phi_n(0,\cdot)}_{L^2} = 1 $ and $ \phi_n(0,\cdot) = 0 $ such that
$$
\int_0^{\infty} | \phi_n (t, te_1) |^2 \dd t \to \infty \qquad \text{as } n \to \infty. 
$$
See also \cite[Prop. 1.1]{Taoexp} for an explicit example of such initial velocity for which the solution has $ \phi(t,t e_1) = \infty $ for all $ t \in (1,2) $ (presented as an adaptation of Stein's counterexample in regards to the spherical maximal function). 

Moreover, \cite{montgomery1998time} shows that, when $ n =3 $, \eqref{W} fails even with $ L^{\infty}_x $ replaced by $ BMO $. This is in contrast to the case $ n \geq 4 $ for which that variant can be obtained from the Keel-Tao endpoint \cite{keel1998endpoint} and the embedding $ \dot{W}^{\frac{n}{r},r} \subset BMO $.
\end{remark}

The proof of Theorem \ref{ThmW} proceeds instead in Fourier space and must use rays of velocity $ >1 $ instead of lightrays.

\subsection{Notation} In what follows $ \hat{f} $ denotes the Fourier transform of $ f $. Operators $ m(D) $, such as $ e^{ i t \vm{D}}, e^{i t \Delta} $, are defined by $ \widehat{m(D) f} (\xi)  = m(\xi) \hat{f} (\xi ) $. Here $ \Delta = - \vm{D}^2 $. Homogeneous Sobolev spaces are defined by $ \dot{H}^s = \vm{D}^{-s} L^2  (\mathbb{R}^n) $ for $ s<\frac{n}{2} $. Space-time norms are defined by 
$$
\vn{F}_{L^q_t L^r_x}^q =  \int_{-\infty}^{\infty}  \vn{F(t,\cdot)}_{L^r_x (\mb{R}^n)} ^q \dd t.
$$
Littlewood-Paley projections $ P_{\leq k} $ are defined by $ \widehat{P_{\leq k} h} (\xi) =  \chi( \frac{ \vm{\xi}}{2^k}) \hat{h}(\xi) $, where $ \chi \in C_c^{\infty}(\mb{R}) $, $ \chi \geq 0 $ and $ \chi( \eta) = 1$ for $ \eta $ near $ 0 $. 
For example, $ \widehat{P_{\leq k} \delta_0} (\xi) =  \chi( \frac{ \vm{\xi}}{2^k})$.  Here $ \delta_0 $ denotes the Dirac delta measure.
$ S^{n-1} \subset \mb{R}^n $ denotes the unit sphere. 
The notation $ X \ls Y $ means $ X \leq C Y $ for an absolute constant $ 0<C < \infty $. When $ Y \ls X \ls Y $ we denote $ X \simeq Y $. 

\section{Proof of Theorem \ref{ThmW}}

 Assume that \eqref{W} is true. Then, by duality, one also has
\be \label{StrD} \tag{W'}
\norm{ \int_{-\infty}^{\infty} e^{-it\vm{D}} \frac{1}{\vm{D}^{\xx}} F(t) \dd t}_{L^2_x} \ls \vn{F}_{L^2_t L^1_x}
\ee
for all $ F \in L^2_t L^1_x \cap L^1_t \dot{H}^{-\xx}_x $.

Take any $ f \in L^1_t \cap L^2_t $ with $ \vn{f}_{L^2}=1 $ such that  $ \int_0^{\infty} | \hat{f} (\eta) |^2 \dd \eta \neq 0 $. Consider also the function $ h \in  L^1_x \cap \dot{H}^{-\xx}_x $ with $ \vn{h}_{L^1_x}=1 $ which remains to be defined. Fix a direction $ \omega \in S^{n-1} $ and a constant $ c>0 $.
Applying  \eqref{StrD} with 
$$ F(t,x) \defeq f(t) h(x-ct \omega) $$  
and using Plancherel's theorem one obtains
$$
\norm{ \frac{1}{\vm{\xi}^{\xx}} \hat{h}(\xi) \int_{-\infty}^{\infty} e^{-it\vm{\xi}} e^{-ict \omega \cdot \xi} f(t) \dd t}_{L^2_{\xi}} \ls 1 
$$
which implies 
$$
\norm{ {\vm{\xi}^{-\xx}} \hat{h}(\xi) \hat{f}(\vm{\xi}+c \omega \cdot \xi)  }_{L^2_{\xi}} \ls 1 .
$$
Consider $ h=\delta_{0}, \ \hat{h} \equiv 1 $, or rather, since $ \delta_0 \notin L^1 \cap \dot{H}^{-\xx} $,  take an approximation to the identity sequence $ h_k=P_{\leq k} \delta_0 $ and pass to the limit. Here $ \delta_0 $ denotes the Dirac delta measure. By the Monotone Convergence theorem one obtains
$$
\int_{\mb{R}^n} \frac{1}{\vm{\xi}^{n-1}} | \hat{f}(\vm{\xi}+c \omega \cdot \xi) |^2 \dd \xi \ls 1 .
$$
Write this integral in spherical coordinates, where 
$$ \vm{\xi}=\lmd, \quad \omega \cdot \xi =\lmd \cos \theta, \quad \theta \in [0,\pi]. $$
This implies 
$$
\int_0^{\infty} \int_0^{\pi} | \hat{f}(\lmd(1+c \cos \theta) ) |^2 (\sin \theta)^{n-2} \dd \theta \dd \lmd \ls 1 .
$$
\begin{remark}
For $ n=3 $ one may take $ c=1 $ and change coordinates to $ u=1+\cos \theta$ to obtain
$$ \int_0^{\infty} \int_0^{2} | \hat{f}(\lmd u ) |^2  \dd u \dd \lmd=\int_0^{2} \frac{1}{u} \dd u \int_0^{\infty} | \hat{f} (\eta) |^2 \dd \eta  \ls 1
$$ 
which gives a contradiction. One may view this as a Fourier version of the physical space construction in Remark \ref{Wn3} above, due to \cite{klainerman1993space}. $ \Delta $
\end{remark}

For $ n \geq 4 $, choosing $ c=1 $ does not work. Instead, take $ c=2 $ and restrict the integral to $ \theta \in [ \frac{ \pi}{2}, \frac{2 \pi}{3}] $. On this interval $ \sin \theta \simeq 1 $. As before, one obtains
$$ 
\int_0^{\infty} | \hat{f} (\eta) |^2 \dd \eta \int_{ \frac{ \pi}{2}}^{ \frac{2 \pi}{3}} \frac{(\sin \theta)^{n-2}}{1+2 \cos \theta }   \dd \theta \ls 1. 
$$
The last integral is logarithmically divergent, being $ \simeq \int_0^{1} \frac{1}{y} \dd y $, thus obtaining a contradiction again, which concludes the proof.
$ \Box $

\bibliographystyle{amsplain}
\nocite{*}
\bibliography{biblioStr}

\providecommand{\bysame}{\leavevmode\hbox to3em{\hrulefill}\thinspace}
\providecommand{\MR}{\relax\ifhmode\unskip\space\fi MR }
\providecommand{\MRhref}[2]{%
  \href{http://www.ams.org/mathscinet-getitem?mr=#1}{#2}
}
\providecommand{\href}[2]{#2}
\begin{thebibliography}{10}

\bibitem{brenner1975p}
Philip Brenner, \emph{On lp- lp' estimates for the wave-equation},
  Mathematische Zeitschrift \textbf{145} (1975), no.~3, 251--254.

\bibitem{fang2006some}
Daoyuan Fang and Chengbo Wang, \emph{Some remarks on strichartz estimates for
  homogeneous wave equation}, Nonlinear Analysis: Theory, Methods \&
  Applications \textbf{65} (2006), no.~3, 697--706.

\bibitem{foschi2005inhomogeneous}
Damiano Foschi, \emph{Inhomogeneous strichartz estimates}, Journal of
  Hyperbolic Differential Equations \textbf{2} (2005), no.~01, 1--24.

\bibitem{ginibre1992smoothing}
Jean Ginibre and Giorgio Velo, \emph{Smoothing properties and retarded
  estimates for some dispersive evolution equations}, Communications in
  mathematical physics \textbf{144} (1992), no.~1, 163--188.

\bibitem{ginibre1995generalized}
\bysame, \emph{Generalized strichartz inequalities for the wave equation},
  Journal of functional analysis \textbf{133} (1995), no.~1, 50--68.

\bibitem{guo2018boundary}
Zihua Guo, Ji~Li, Kenji Nakanishi, and Lixin Yan, \emph{On the boundary
  strichartz estimates for wave and schr{\"o}dinger equations}, Journal of
  Differential Equations \textbf{265} (2018), no.~11, 5656--5675.

\bibitem{kapitanski1989some}
Lev~Vil'evich Kapitanski, \emph{Some generalizations of the strichartz--brenner
  inequality}, Algebra i Analiz \textbf{1} (1989), no.~3, 127--159.

\bibitem{keel1998endpoint}
Markus Keel and Terence Tao, \emph{Endpoint strichartz estimates}, American
  Journal of Mathematics \textbf{120} (1998), no.~5, 955--980.

\bibitem{klainerman1993space}
S~Klainerman and M~Machedon, \emph{Space-time estimates for null forms and the
  local existence theorem}, Communications on Pure and Applied Mathematics
  \textbf{46} (1993), no.~9, 1221--1268.

\bibitem{klainermanlecture}
Sergiu Klainerman, \emph{Lecture notes in analysis (2011)},
  https://web.math.princeton.edu/~seri/homepage/courses/Analysis2011.pdf.

\bibitem{klainerman1998algebraic}
Sergiu Klainerman and Matei Machedon, \emph{On the algebraic properties of the
  $ h_{n/2, 1/2} $ spaces}, International Mathematics Research Notices
  \textbf{1998} (1998), no.~15, 765--774.

\bibitem{lindblad1995existence}
Hans Lindblad and Christopher~D Sogge, \emph{On existence and scattering with
  minimal regularity for semilinear wave equations}, Journal of Functional
  Analysis \textbf{130} (1995), no.~2, 357--426.

\bibitem{machihara2005endpoint}
Shuji Machihara, Makoto Nakamura, Kenji Nakanishi, and Tohru Ozawa,
  \emph{Endpoint strichartz estimates and global solutions for the nonlinear
  dirac equation}, Journal of Functional Analysis \textbf{219} (2005), no.~1,
  1--20.

\bibitem{mockenhaupt1993local}
Gerd Mockenhaupt, Andreas Seeger, and Christopher~D Sogge, \emph{Local
  smoothing of fourier integral operators and carleson-sjolin estimates},
  Journal of the American Mathematical Society (1993), 65--130.

\bibitem{montgomery1998time}
SJ~Montgomery-Smith, \emph{Time decay for the bounded mean oscillation of
  solutions of the schr{\"o}dinger and wave equations}, Duke Mathematical
  Journal \textbf{91} (1998), no.~2, 393--408.

\bibitem{segal1976space}
Irving Segal, \emph{Space-time decay for solutions of wave equations}, Advances
  in Mathematics \textbf{22} (1976), no.~3, 305--311.

\bibitem{sogge1995lectures}
Christopher~Donald Sogge, \emph{Lectures on non-linear wave equations}, vol.~2,
  International Press Boston, MA, 1995.

\bibitem{sterbenz2005angular}
Jacob Sterbenz, \emph{Angular regularity and strichartz estimates for the wave
  equation}, International Mathematics Research Notices \textbf{2005} (2005),
  no.~4, 187--231.

\bibitem{strichartz1970priori}
Robert~S Strichartz, \emph{A priori estimates for the wave equation and some
  applications}, Journal of Functional Analysis \textbf{5} (1970), no.~2,
  218--235.

\bibitem{strichartz1977restrictions}
Robert~S Strichartz et~al., \emph{Restrictions of fourier transforms to
  quadratic surfaces and decay of solutions of wave equations}, Duke
  Mathematical Journal \textbf{44} (1977), no.~3, 705--714.

\bibitem{Taoexp}
Terence Tao, \emph{Counterexamples to the n=3 endpoint strichartz estimate for
  the wave equation}, Expository note:
  https://www.math.ucla.edu/~tao/preprints/Expository/stein.dvi.

\bibitem{tao2006counterexample}
\bysame, \emph{A counterexample to an endpoint bilinear strichartz inequality},
  Electronic Journal of Differential Equations \textbf{2006} (2006), no.~151,
  1--6.

\bibitem{tao2006nonlinear}
\bysame, \emph{Nonlinear dispersive equations: local and global analysis}, no.
  106, American Mathematical Soc., 2006.

\bibitem{yajima1987existence}
Kenji Yajima, \emph{Existence of solutions for schr{\"o}dinger evolution
  equations}, Communications in Mathematical Physics \textbf{110} (1987),
  no.~3, 415--426.

\end{thebibliography}

\end{document}